\documentclass[11pt]{amsart}
\usepackage{amsthm}
\usepackage{amssymb}
\usepackage{amsmath}

\newtheorem{theorema}{Theorem}
\newtheorem{theorem}{Theorem}[section]
\newtheorem*{theorem*}{Theorem}
\newtheorem{lemma}[theorem]{Lemma}

\newtheorem{proposition}[theorem]{Proposition}
\newtheorem*{proposition*}{Proposition}

\theoremstyle{remark}
\newtheorem{remark}[theorem]{Remark}
\newtheorem*{remark*}{Remark}

\numberwithin{equation}{section}
\newcommand{\st}[1]{\ensuremath{^{\scriptstyle \textrm{#1}}}}

\gdef\myletter{}

\let\savetheequation\theequation
\def\theequation{\savetheequation\myletter}

\def    \Cinf   {C^\infty}

\newcommand{\CC}{{\mathbb C}}
\newcommand{\C}{{\mathbb C}}  

\newcommand{\RR}{{\mathbb R}}
\newcommand{\R}{{\mathbb R}}  
\newcommand{\ZZ}{{\mathbb Z}}
\newcommand{\bbT}{{\mathbb T}} 

\newcommand{\fg}{{\mathfrak g}}
\newcommand{\fh}{{\mathfrak h}} 

\newcommand{\A}{{\mathcal A}}
\newcommand{\B}{{\mathcal B}}

\newcommand{\cut}{\operatorname{cut}}
\newcommand{\even}{\operatorname{even}}  
\newcommand{\Hom}{\operatorname{Hom}}    
\newcommand{\red}{\operatorname{red}} 
\newcommand{\res}{\operatorname{res}}
\renewcommand{\span}{\operatorname{span}}  
\newcommand{\vol}{\operatorname{vol}}
\newcommand{\SL}{\operatorname{SL}}

\def    \inv    {^{-1}}

\newcommand{\alphaparenlist}{
  \renewcommand{\theenumi}{\alph{enumi}}%
  \renewcommand{\labelenumi}{(\theenumi)}%
}

\begin{document}

\title[Melrose--Uhlmann projectors]{Melrose--Uhlmann projectors, the
  metaplectic representation and symplectic cuts}

\author{V. Guillemin and E. Lerman}\thanks{Supported in
part by NSF grants  DMS-0104116 (VG) and DMS-980305 and
DMS-0204448 (EL)} 


\maketitle

\section*{Abstract}

By applying the symplectic cutting operation to cotangent
bundles, one can construct a large number of interesting
symplectic cones.  In this paper we show how to attach algebras
of pseudodifferential operators to such cones and describe the
symbolic properties of the algebras.

\section*{Introduction}

The Melrose--Uhlmann projectors which we refer to in the title of
this article are projection operators which look microlocally
like the standard Szeg\"o projectors on $L^2(S^1)$.  They belong
to a class of pseudodifferential operators with singular symbols
which were studied by Melrose--Uhlmann in \cite{MU} and by one of
us in \cite{Gu}.  One of the main goals of this paper will be to
give a microlocal description of the algebra of classical
pseudodifferential operators which commute with such a projection
operator.

Another of the main goals of this paper will be to examine some
microlocal aspects of a basic operation in cobordism theory:  the
\emph{cutting} operation.  Let $M$ be a $\Cinf$ manifold, $\tau :
S^1 \times M \to M$ an action of $S^1$ on $M$ and $\Phi : M \to
\RR$ an $S^1$-invariant function.  If zero is a regular value of
$\Phi$ the set
\begin{displaymath}
  W= \{ p \in M \, , \,\, \Phi (p) \geq 0 \}
\end{displaymath}
is a manifold with boundary, and if $S^1$ acts freely on the
boundary, one gets a $\Cinf$ manifold without boundary by
collapsing the circle orbits in the boundary to points.  This new
manifold, which we will denote by $M$, is the disjoint union of
the manifold, $M_{\red} = \Phi^{-1} (0) /S^1$ and the interior,
$W^0$, of $W$; and $M_{\red}$ sits inside $M$ as a codimension
$2$ submanifold.  For example let $M=\CC^n$ and let $\tau$ be
multiplication by unit complex numbers.  If $\Phi (z) =
|z|^2-1$, then $M$ is the blow up of $\CC^n$ at $0$ and $M_{\red}
= \CC P^{n-1}$ is the exceptional divisor.  On the other hand if
$\Phi (z) = 1-|z|^2$, then $M$ is $\CC P^n$ and $M_{\red} = \CC
P^{n-1}$.

It was observed several years ago by one of us (see \cite{Le})
that this cutting operation can be symplecticized.  Namely
suppose that $M= (M,\omega)$ is a symplectic manifold, $\tau$ a
Hamiltonian action and $\Phi$ the moment map associated with this
action.  Then the symplectic form on $W^0=M-M_{\text{red}}$ extends
smoothly to a symplectic form on $M_{\text{cut}}$ and so also does the
action $\tau$ and moment map, $\Phi$.  Moreover, $M_{\text{red}}$ is a
symplectic submanifold of $M_{\text{cut}}$ and, as an abstract
symplectic manifold, is  isomorphic to the usual
symplectic  reduction of $M$ by $\tau$.

To prove these assertions one needs a somewhat different
description of $M_{\text{cut}}$.  Consider the product manifold, $M
\times \CC$, with the product symplectic form, $\omega_M - \omega _{\CC}$,
and the action on it of $S^1 \times S^1$.  The moment map for
this product action is $(\Phi (m), - |z|^2)$; so if we restrict
to the diagonal subgroup of $S^1 \times S^1$ we get a Hamiltonian
action of $S^1$ on $M \times \CC$ with moment map, $\Psi (m,z) =
\Phi (m) - |z|^2$, and it is not hard to see that $M_{\text{cut}}$ can be
identified with the reduced space
\begin{equation}
  \label{eq:0.3} (M \times \CC)_{\text{red}} = \Psi^{-1} (0) /S^1 \, .
\end{equation}
Moreover this space has a residual action on it of $S^1$, and it
is not hard to see that this action coincides with the action of
$S^1$ described above.

Suppose now that the action $\tau$ can be quantized;
i.e., suppose that one can associate with $(M,\tau)$ a
representation, $\tau^{\#}$, of $S^1$ on a Hilbert space, $Q(M)$,
by some kind of ``quantization'' procedure.  Then, in view of the
fact that the symplectic form on $M \times \CC$ defined above is
the product of the symplectic form on $M$ and on  $\CC$, one gets for the
quantization of $M \times \CC$
\begin{displaymath}
  Q(M) \otimes Q (\CC)^*
\end{displaymath}
or equivalently
\begin{displaymath}
  \Hom (Q (\CC) , \, Q(M)) \, .
\end{displaymath}
Thus by the ``quantization commutes with reduction'' principle
one gets for the reduced space, $M_{\cut} = (M \times
\CC)_{\red}$ the quantization
\begin{displaymath}
 \Hom (Q (\CC) ,\, Q(M))^{S^1}  \, .
\end{displaymath}
To complete this quantum description of $M_{\cut}$ we still have
to specify a quantization $Q(\CC)$, of the action of $S^1$ on
$\CC$ and for this there is a more or less canonical candidate,
the oscillator representation of $S^1$ on $L^2 (\RR)$.  Thus the
Hilbert space
\begin{equation}
  \label{eq:0.4}
  \Hom (L^2 (\RR) \, , Q(M))^{S^1}
\end{equation}
is an obvious candidate for $Q(M_{\cut})$.\footnote{To make
  \ref{eq:0.4} into a Hilbert space we will take the intertwining
  operators in this ``$\Hom$'' to be Hilbert--Schmidt.}

To see how this construction is related to the theory of
Melrose--Uhlmann projectors let $H_n \, , \, n=0,1, \ldots$ be
the one-dimensional subspace of $L^2 (\RR)$ spanned by the
$n$\st{th} Hermite function, $h_n$.  This subspace transforms as
$e^{in\theta}$ under the action of $\theta \in S^1$.  Therefore
the space
\begin{displaymath}
  \Hom (H_n  , \, Q(M))^{S^1}
\end{displaymath}
can be identified with the space
\begin{displaymath}
  Q_n (M) = \{ f \in Q(M) \, , \,
  \tau^{\#} (e^{i\theta}) f = e^{in\theta} f \}
\end{displaymath}
via the map
\begin{displaymath}
  T \mapsto Th_n \, ,
\end{displaymath}
and the space (\ref{eq:0.4}) can be identified with the direct
sum
\begin{equation}
  \label{eq:0.5}
\bigoplus^{\infty}_{n=0} Q_n (M) \, .
\end{equation}
Let us denote by $\Pi_+$ the orthogonal projection of $Q(M)$ onto
the space (\ref{eq:0.5}).  The examples we will be interested in
in this paper with be quantizations defined using microlocal
analysis, and for these examples $\Pi_+$ will be a projector of
Melrose--Uhlmann type.  Moreover in these examples there will be
a natural algebra of ``quantum observables'' on $M$:  either
pseudodifferential operators or Toeplitz operators, and hence a
natural algebra of quantum observables on $M_{\cut}$, namely the
operators which commute with $\Pi_+$.

Finally we'll explain why the metaplectic representation is
involved in the construction we've just described.  The
oscillator representations of $S^1$ on $L^2 (\RR)$ is
unfortunately not a representation of $S^1$ itself but of its
metaplectic double cover.  This double cover is just another copy
of $S^1$; so there would seem to be no problem in substituting it
for $S^1$ in the definition (\ref{eq:0.4}).  However, if one
wants to attach symbols to the quantum observables we just
defined, the fact that the $S^1$ acting on $\CC$ is not the same
$S^1$ as that acting on $L^2 (\RR)$ causes some unpleasant parity
complications and one has to make use of metaplectic techniques
to deal with these complications.

A few words about the contents of this article.  For simplicity
we will henceforth assume that the manifold $M$ above is the
cotangent bundle of a compact manifold, $X$, and that the algebra
of ``quantum observables'' is the algebra
of pseudodifferential operators, $\Psi (X)$.\footnote{However,
    most of the results below are true, mutatis mutandis, for the
    algebra of Toeplitz operators on a strictly pseudoconvex domain.}
As for the  action, $\tau$ we will assume it is a \emph{canonical} action,
  i.e.,~each of the symplectomorphisms, $\tau (e^{i\theta})$, is
  a canonical transformation.  By a theorem of de la
  Harpe--Karoubi \cite{HK} every such action can be quantized by
  a unitary representation
  \begin{displaymath}
    \tau^{\#} : S^1 \to U (H) \, , \, H=L^2 (X) \, ,
  \end{displaymath}
by Fourier integral operators; and for this representation the
projector $\Pi_+$ is of Melrose--Uhlmann type.
(See \cite{Gu} theorem~4.4.  We will also prove this explicitly in
\S4 by showing that $\Pi_+$ is microlocally conjugate to the
standard Szeg\"o projector.)
The main result of this article is the following
\begin{quote}
``{\bf Theorem }''
  Let $\Psi_+$ be the algebra of pseudodifferential operators
  which commute with $\Pi_+$.  Then the algebra $\Pi_+ \Psi_+
  \Pi_+$ quantizes the algebra of classical observables, $\Cinf
  (M_{\text{cut}})$.
\end{quote}
The second statement needs some amplification (which will be
supplied in \S5); however the reason for the quotation marks is
the parity complications we referred to above.  We will
discuss this ``metaplectic glitch'' in more detail in \S1 and  will show that there are two
ways of dealing with it:  one by making the action of $S^1$ on
$M$ a ``metaplectic'' action and the other by making the action
of $S^1$ on $\C^1$ a ``metaplectic'' action.  We will show
that both these alternatives give rise to an interesting symbol
calculus for operators in $\Pi_+ \Psi_+ \Pi_+$.

In section~2 we will discuss a differential operator version of
the ``theorem'' above for the manifold $X=S^1$ and the standard
Szeg\"o projector, and then in section~3 we will extend this result
to the algebra of pseudodifferential operators on product
manifolds of the form, $X=Y \times S^1$.  In section~4 we will
show that it suffices to prove our ``theorem'' in this case by
showing that there exists a Fourier integral operator locally
conjugating the general case to this case.  Finally in section~5
we will discuss the symbolic calculus of the algebra $\Pi \Psi_+
\Pi$.  We will show that an operator of degree $r$ in this
algebra has a leading symbol which is an homogeneous function of
degree $r$ on $M_{\text{cut}}$ and that products and Poisson brackets of
symbols correspond to products and commutators of operators.  We
will also show that this algebra can be equipped with a residue
trace which, for operators of degree $-d$, $d=\dim M/2$, is given
by integrating the leading symbol of the operator over $M_+$, and
will deduce from this a Weyl law for operators of elliptic type.

Finally in section~\ref{sec:6} we will discuss what happens when one
starts with a cotangent bundle and applies to it repeated symplectic
cuts. One can construct in this way a lot of interesting symplectic
cones, and by the techniques of this paper one gets (modulo the
$\ZZ_2$ problems discussed above) algebras of classical
polyhomogeneous pseudodifferential operators quantizing these cones.
The details of this construction will be spelled out elsewhere but in
section~\ref{sec:6} we will indicate (roughly) how to quantize in this
way the cones over the classical three dimensional lens spaces.

\section{The metaplectic glitch}
\label{sec:intro}

Let $M$ be a manifold with an action
$\tau$ of a circle $S^1$ and an $S^1$ invariant function $\Phi: M\to
\RR$.  Suppose $S^1$ acts freely on the level set $\Phi\inv (0)$.  Then
the quotient $M_{\red}: = \Phi\inv (0)/S^1$ is a manifold.  Consider the
manifold with boundary $\{m\in M \mid \Phi (m) \geq 0\}$, and collapse
the circle orbits in the boundary to points.  The resulting space
$$
M_{\cut} := \{m\in M \mid \Phi (m) \geq  0\} /\!\!\sim \, ,
$$
where $\sim$ is the relation described above (cf.\ (\ref{eq.rel})
below), is a $C^0$ manifold.  The manifold $M_{\red}$ embeds naturally
in $M_{\cut}$ as a codimension 2 submanifold and the difference
$M_{\cut} \smallsetminus M_{\red}$ is homeomorphic to $\{m\in M \mid
\Phi (m) > 0\}$.

If, in addition, $M$ is a symplectic manifold, the action $\tau$ is
Hamiltonian and $\Phi :M\to \RR$ is the corresponding moment map then
$M_{\cut}$ is symplectic.  More specifically

\begin{proposition}\label{prop:cut}
Let $(M, \omega)$ be a symplectic manifold with a Hamiltonian action
$\tau$ of $S^1$; let $\Phi : M \to \RR$ denote a corresponding moment
map.  Suppose $S^1$ acts freely on $\Phi\inv (0)$.  Define an
equivalence relation $\sim$ on $\{m \in M \mid \Phi (m) \geq 0\}$ for
$m\not = m'$ by the identification:
\begin{equation} \label{eq.rel}
m\sim m' \Longleftrightarrow \Phi (m) = \Phi (m') = 0 \,\,
\mbox{ and } \,\, m = \lambda \cdot m' \,\,
\mbox{ for some } \,\, \lambda \in S^1 .
\end{equation}
Then
\begin{enumerate}
\item The $C^0$ manifold $M_{\cut}$ can be given the structure of a $C^\infty$
symplectic manifold $(M_+, \omega_+)$ so that the reduced space
$M_{\red} = \Phi\inv (0)/S^1$ embeds symplecticly and the difference $M_+
\smallsetminus M_{\red}$ is symplectomorphic to
$\{m \in M \mid \Phi (m) > 0\}$.
\item  Alternatively, the $C^0$ manifold $M_{\cut}$ can be given the structure
of a $C^\infty$ symplectic orbifold $(M_{++}, \omega_{++})$ so that
the set of regular points is symplectomorphic to $\{m \in M \mid \Phi
(m) > 0\}$, the set of singular points is symplectomorphic to the
reduced space $M_{\red}$, and the structure group
of all points in $M_{\red}$ is $\ZZ_2$.
\end{enumerate}
\end{proposition}
\begin{remark}
Even though $M_+$ and $M_{++}$ are the same as topological spaces,
namely $M_{\cut}$, they are not the same as orbifolds.  In particular
$\Cinf (M_+) \neq
\Cinf (M_{++})$.
\end{remark}
\begin{remark}
One readily sees from the proof below that the Hamiltonian action
$\tau$ of $S^1$ on $(M, \omega)$ descends to a Hamiltonian action of
$S^1$ on $(M_+, \omega_+)$ which fixes $M_{\red}$ pointwise and makes
the embedding $\{\Phi >0\} \hookrightarrow M_+$ equivariant.  The same
statement holds for $M_{++}$.
\end{remark}

\begin{proof}
Consider the diagonal action of $S^1$ on $(M\times \CC, \omega -
idz\wedge d\bar{z})$.  The map $\tilde{\Phi} (m, z) = \Phi (m) - |z|^2$
is a corresponding moment map.  Since $S^1$ acts freely on $\Phi\inv
(0)$ it acts freely on $\tilde{\Phi} \inv (0)$.  Hence $M_+ : =
\tilde{\Phi} \inv (0)/S^1$ is a symplectic manifold.  The composition
of the embedding $j: \{\Phi \geq 0\} \hookrightarrow \tilde{\Phi} \inv
(0)$, $j (m) = (m, \sqrt{\Phi (m)})$ with the orbit map $\tilde{\Phi}
\inv (0) \to \tilde{\Phi}\inv (0)/S^1 = M_+$ is onto.  It induces a
homeomorphism $\varphi: M_{\cut} = \{\Phi \geq 0\}/\!\!\sim \,\, \to M_+$.
Note that $\varphi|_{\{\Phi > 0\}}$ is an open embedding.  Moreover,
since $j^* (\omega -idz\wedge d\bar{z}) =
\omega$, it is symplectic.  Similarly one checks that the difference
$M_+ \smallsetminus \varphi (\{\Phi > 0\}) $ is the reduced space
$M_{\red}$.  This proves the first part of the theorem.

Denote elements of $\CC /\ZZ_2$ by $[z]$, so that $[z]=[-z]$.
Consider the $S^1$ action on $\CC/\ZZ_2$ given by $\mu \cdot [z] =
[\sqrt{\mu} z]$. This action is well-defined and preserves the
symplectic form $-[i \, dz \wedge d\bar{z}]$ on $\CC /\ZZ_2$
corresponding to $-i \, dz \wedge d\bar{z}$; the moment map for
this action is the map, $[z]\mapsto -|z|^2$.  Now consider the
diagonal action of $S^1$ on $(M\times \CC/\ZZ_2, \omega -[i \, dz
\wedge d\bar{z}])$ and proceed as in the first part of the proof, denoting
the reduction of $M\times \CC/\ZZ_2$ at zero by $M_{++}$.
\end{proof}

Now let $M$ be a cotangent bundle of a compact manifold, $X$, of
dimension $n>1$; and let the action, $\tau$, above be a canonical
action (an action preserving the canonical cotangent one-form, $\Sigma
\xi_i \, dx_i$).  By the theorem of
de~la~Harpe--Karoubi--Weinstein that we cited in the introduction,
there exists a representation, $\tau^\#$, of $S^1$ on $L^2 (X)$ which
quantizes $\tau$ in the sense that for each $e^{i\theta} \in S^1$,
$\tau^\# (e^{i\theta})$ is a unitary Fourier integral operator with
$\tau (e^{i\theta})$ as its underlying canonical transformation.  Let
\begin{displaymath}
  \Pi : L^2 (X) \to L^2 (X)
\end{displaymath}
be orthogonal projection onto the space
\begin{equation}
  \label{eq:1.1}
  \span \{ f \in L^2 (X) \mid \quad
  \tau^\# (e^{i\theta}) f
  = e^{in\theta} f, \quad n \geq 0 \} \, ,
\end{equation}
and let $\Psi_+$ be the algebra of classical pseudodifferential
operators which commute with $\Pi$.  The main result of this
paper (modulo a few qualifications which we will explain shortly)
 asserts:
\begin{gather}
\text{ \em $\Pi$ is a projector of Melrose--Uhlmann type,
and the algebra $\Pi \Psi_+ \Pi$} \tag{$*$}\\
\text{ \em quantizes the algebra of classical observables,
$\Cinf (M_+)$.} \notag
\end{gather}

As we remarked in the introduction
there is a metaplectic glitch involved in making
the statement above
correct.  To explain this glitch we note that, since
$M=T^* X$, the obvious candidate for the quantum Hilbert space
to associate with $M$ is $L^2(X)$; and since
\begin{displaymath}
  \CC = \RR^2 = T^* \RR
\end{displaymath}
the obvious candidate for the quantum Hilbert space to associate
with $\CC$ is $L^2 (\RR)$.  Thus, if one subscribes to the
principle that ``quantization commutes with reduction'' one should
associate with $M_+$ the quantum Hilbert space
\begin{equation}
  \label{eq:1.2}
  \Hom (L^2 (\RR), L^2 (X))^{S^1} \, .
\footnote{There is a ``Hom'' rather than a tensor product here
  because the symplectic cutting procedure requires one to take
  the symplectic form on $\CC$ to be the negative of the usual
  symplectic form  (cf.~proof of  Proposition~\ref{prop:cut} above).}
\end{equation}
We must, of course, still specify how $S^1$ is to act on $L^2
(\RR)$ for (\ref{eq:1.2}) to make sense; and this, we will see,
is the source of the ``metaplectic glitch'' that we referred to
above.  Let's briefly review how the metaplectic (or
Segal--Shale--Weil) representation of $S^1$ on $L^2 (\RR )$ is
defined:  Let $x$ and $y$ be the standard Darboux coordinates on
$\RR^2$ and let $\fh_3 = \span \{ x,y,1 \}$.  This space sits
inside the Poisson algebra, $\Cinf (\RR^2)$, as a
three-dimensional Heisenberg algebra, and can be represented on
$L^2 (\RR)$ by the standard Schroedinger representation
\begin{equation}
  \label{eq:1.3}
  x \to x, \quad y \to \tfrac{\partial}{\partial x} ,
  \quad 1 \to I \, .
\end{equation}
This exponentiates to a representation, $\kappa$, of the Heisenberg
group, $H_3$, on $L^2 (\RR)$; and, by the Stone--Von~Neumann theorem,
$\kappa$ is the unique irreducible representation of $H_3$ for which
the center, $\RR$, of $H_3$ acts as $e^{it}I$.  Consider now the
symplectic action of $S^1$ on $\RR^2$ given by $\theta \to
e^{i\theta}$.  Being a linear action this preserves $\fh_3$, and being
symplectic, acts on $\fh_3$ by Lie algebra automorphisms.  Hence,
since $H_3$ is simply connected, this action can be exponentiated to
an action, $\rho$, of $S^1$ on $H_3$ by Lie group automorphisms; and
this enables one to define, for every $\theta$, a new representation,
$\kappa_{\theta}$, of $H_3$ on $L^2 (X)$ by setting
\begin{displaymath}
  \kappa_{\theta} (h) = \kappa(h_{\theta}), \quad
  h_{\theta} = \rho (e^{i\theta}) h \, .
\end{displaymath}
This representation is identical with $\kappa$ on the center of  $H_3$;
so by the Stone--Von~Neumann theorem $\kappa$ and $\kappa_{\theta}$ are
isomorphic:  there exists a unitary operator
\begin{displaymath}
  \gamma_{\theta} : L^2 (X) \to L^2 (X)
\end{displaymath}
such that $\gamma^{-1}_{\theta} \kappa \gamma_{\theta} = \kappa_{\theta}$.
Moreover, since $\kappa$ is irreducible, this operator is unique up to
a constant multiple of module one.  From this uniqueness it is
easy to see that $\gamma_{\theta_1 + \theta_2}$ is a constant
multiple of $\gamma_{\theta_1} \gamma_{\theta_2}$; i.e.,~the map
\begin{equation}
  \label{eq:1.4}
  e^{i\theta} \to \gamma (\theta)
\end{equation}
is a \emph{projective} representation of $S^1$ on $L^2 (X)$.  The
problem of converting this projective representation into a
bona fide \emph{linear} representation is a standard problem in
representation theory and involves an obstruction which sits in
the group cohomology of the group, $S^1$.  For (\ref{eq:1.4})
this obstruction unfortunately doesn't vanish; but one can make
it vanish by pulling it back to the metaplectic double cover,
$\tilde{S}^1$, of $S^1$.  Since $\tilde{S}^1$ is just the group
$S^1$ itself, double covering itself by the map
\begin{equation}
  \label{eq:1.5}
  e^{i\theta} \mapsto e^{2i \theta}
\end{equation}
one gets a linear representation, $\tilde{\gamma}$, of $S^1$ on
$L^2 (\RR)$ by composing (\ref{eq:1.4}) with (\ref{eq:1.5}) and
adjusting constant multiples.  This is, by definition, the
metaplectic representation of $S^1$ on $L^2 (\RR)$; and its
clear from this definition that its the \emph{only}
representation of $S^1$ on $L^2 (\RR)$ compatible with (\ref{eq:1.3}).

Coming back to the space of intertwining operators (\ref{eq:1.2}), if
the representation of $S^1$ on $L^2 (\RR)$ is the metaplectic
representation, the space (\ref{eq:1.2}) is not strictly speaking
well-defined since the ``$S^1$'' acting on $L^2 (\RR)$ is \emph{not}
the same group as the ``$S^1$'' acting on $L^2 (X)$ and on $M \times
\CC$: it is the metaplectic double cover of this group.  This is the
``metaplectic glitch'' which we referred to above.  We will discuss
below two ways of dealing with it, one of which leads to an
interesting quantization of $M_+$ and the other to an interesting
quantization of $M_{++}$.

The first way is to make the action of $S^1$ on the second factor
of (\ref{eq:1.2}) a metaplectic action.  Namely let $\ZZ_2 = \{
\pm 1 \} = \{ \lambda \in S^1 , \lambda^2 = 1 \}$.  Then $S^1
/\ZZ_2$ acts on $M/\ZZ_2$, and the quantization of this action is
the action of $S^1 / \ZZ_2$ on $L^2 (X/{\ZZ_2})= L^2 (X)^{\ZZ_2}$.  Let's
temporarily relabel the groups, $S^1$ and $S^1 /\ZZ_2$, letting
$S^1$ temporarily be labeled $\tilde{S}^1$ and $S^1/\ZZ_2$,
temporarily labeled $S^1$; and let's replace the space of
intertwining operators, (\ref{eq:1.2}), by
\begin{equation}
  \label{eq:1.6}
  \Hom (L^2 (\RR), L^2 (X)^{\ZZ_2})^{\tilde{S}^1} \, ,
\end{equation}
which is now well-defined since the same group is acting on both
factors.  Let $h_i \in L^2 (\RR )$ be the $i$\st{th} Hermite
function, normalized to have $L^2$-norm one.  Then if $T$ is an
intertwining operator belonging to the space~(\ref{eq:1.6}),
$Th_n=0$ for all $n$ odd and the map
\begin{equation}
\label{eq:1.7}
           T \to \sum^{\infty}_{n=0} Th_{2n}
\end{equation}
maps the space (\ref{eq:1.6}) bijectively onto the
$\ZZ_2$-invariant part of the space (\ref{eq:1.1}).  Let
$\Pi^{\even}$ be the orthogonal projection of $L^2 (X)^{\ZZ_2}$
onto this space and let $\Psi^{\even}_+$ be the
ring of $\ZZ_2$-invariant classical pseudodifferential operators
which commute with $\Pi^{\even}$.  Then the following even
version of assertion ($*$) above is true:

\begin{theorema}
  \label{th:1}
$\Pi^{\even}$ is a projector of Melrose--Uhlmann type and the algebra
$\Pi^{\even} \Psi^{\even}_+ \Pi^{\even}$ quantizes the algebra of
classical observables, $\Cinf (M_+)_{\even}$.
\end{theorema}

The second way of dealing with this ``metaplectic glitch'' is to
make the action of $\tilde{S}^1$ on the first factor of
(\ref{eq:1.2}) an action of $S^1$ by noting that one gets from
the metaplectic representation a representation of $\tilde{S}^1
/\ZZ_2 =S^1$ on $L^2 (\RR)^{\ZZ_2}$.  This makes the space of
intertwining operators
\begin{equation}
  \label{eq:1.8}
  \Hom (L^2 (\RR )^{\ZZ_2},  L^2 (X))^{S^1}
\end{equation}
well-defined.  Moreover, the action of $S^1$ on $L^2
(\RR)^{\ZZ_2}$ is given by $e^{i\theta} \bullet h_{2n}=
e^{in\theta}h_{2n}$, so the mapping (\ref{eq:1.7}) maps the space
(\ref{eq:1.8}) bijectively onto the space (\ref{eq:1.1}), and the
projector, $\Pi$, is projection onto its image.  Now, however,
the classical counterpart of the space (\ref{eq:1.8}) is no
longer $M_+$ but $M_{++}$. Indeed, the first factor in (\ref{eq:1.8}) is
the space, $L^2 (\RR)^{\ZZ_2}$, which one can think of as the
quantization of the orbifold, $(T^*\RR)/\ZZ_2$.  Therefore, by
the principle of ``quantization commutes with reduction'' the classical
counterpart of (\ref{eq:1.8}) is the symplectic reduction at zero
of the orbifold, $M \times ((T^*\RR)/\ZZ_2)$, i.e., it is $M_{++}$.  Our second
version of the ``theorem'' above states

\begin{theorema}\label{th:2}
  $\Pi \Psi_+ \Pi$ is the quantization of the
  algebra of classical observables, $\Cinf (M_{++})$.
\end{theorema}

We now sketch the proof of Theorems~\ref{th:1} and \ref{th:2} for the
space $M=T^*S^1$ (with pseudodifferential operators replaced by
differential operators).  By definition $(T^*S^1)_+$ is the reduction
at zero of the manifold $(T^* S^1 \times \CC = \RR \times S^1 \times
\CC, \, ds \wedge \frac{d\lambda}{i\lambda} - i dz \wedge d\bar{z})$
(where $(s, \lambda
= e^{i\theta}, z)\in \RR \times S^1 \times \CC$), by the $S^1$ action
$$
\mu \cdot (s, \lambda, z) = (s, \mu \lambda, \mu z)
$$
with moment map
$$
\tilde{\Phi} (s, \lambda, z) = s - |z|^2.
$$
The set $\{ (|z|^2, 1, z) \in \RR \times S^1 \times \CC \mid z\in
\CC\}$ parameterizes $S^1$ orbits in $\Psi\inv (0)$.  Hence the map
$\pi :\tilde{\Phi}\inv (0) \to \CC$, $\pi (s = |z|^2, \lambda, z) =
\lambda\inv z$ induces a diffeomorphism $\tilde{\Phi}\inv (0)/S^1 \to \CC$.
The embedding $j: [0,\infty) \times S^1 \to \tilde{\Phi} \inv (0)$, $j(s,
\lambda) = (s, \lambda, \sqrt{s})$ has the property that the
composition $\sigma = \pi\circ j: [0, \infty ) \times S^1 \to \CC$ is
onto; it is one-to-one on $(0,\infty) \times S^1$ and maps $\{0\}
\times S^1$ to $0$. Note that $\sigma (s, \lambda ) = \lambda\inv \sqrt{s}$
and that $\sigma$ induces a homeomorphism $\varphi: \left(([0,\infty) \times
S^1)/\sim \right)\to \CC$.

Now consider the ring of real $\CC$-valued polynomials on $\RR^2 =
\CC$ invariant under the $\ZZ_2$ action $z\mapsto -z$.  It is
generated by $z^2, |z|^2$ and $\bar{z}^2$.  Note that $\sigma^* z^2 =
\lambda ^{-2} s = e^{-2i\theta} s$, $\sigma ^* |z|^2 = s$ and $\sigma
^* \bar{z}^2 = e^{2i\theta} s$.  On the other hand we will show in
\S~\ref{sec:2} that the ring of differential operators on $S^1$ which
commute with the projector, $\Pi^{\even}$, is generated by the
operators
\begin{displaymath}
  \frac{1}{i} \frac{d}{d\theta} e^{2i\theta}, \quad
  \frac{1}{i} e^{-2i\theta} \frac{d}{d\theta} \hbox{ and }
  \frac{1}{i}\frac{d}{d\theta} \, ;
\end{displaymath}
and the symbols of these operators are exactly $se^{2i\theta}$,
$se^{-2i\theta}$ and $s$.  Thus the ring of even polynomial
functions on $(T^*S^1)_+$ is exactly the ring of the symbols of
differential operators which commute with $\Pi^{even}$.

The proof of Theorem~\ref{th:2} is similar.  First note
that
\begin{displaymath}
  (T^*\RR)/\ZZ_2 = \RR^2 /\ZZ_2 = \CC /\ZZ_2 \, .
\end{displaymath}
Let's again denote elements of $\CC /\ZZ_2$ by $[z]$, so that
$[z]=[-z]$.  Consider the $S^1$ action on $\CC/\ZZ_2$ given by $\mu
\cdot [z] = [\sqrt{\mu} z]$. As we noted previously this action is
well-defined and preserves the symplectic form on $\CC /\ZZ_2$
corresponding to $-i \, dz \wedge d\bar{z}$; and the moment map for
this action is the map, $[z]\mapsto -|z|^2$.

Now let's check what $(T^* S^1)_{++}$ looks like.  By definition
$(T^*S^1)_{++}$ is the reduction at zero of the orbifold
\begin{displaymath}
  (T^*S^1) \times (\CC /\ZZ_2) = \RR \times S^1
     \times (\CC /\ZZ_2)
\end{displaymath}
by the circle action
\begin{displaymath}
  \mu \cdot  (s,\lambda = e^{i\theta}, [z])
    = (s, \mu \lambda , [\sqrt{\mu}z]) \, ,
\end{displaymath}
the moment map for this action being the function, $\tilde{\Phi} (s,\lambda
,[z]) =s-|z|^2$.  Arguing as above we get a surjective map $\sigma  :
[0,\infty) \times S^1 \to \CC /\ZZ_2$ which is one-to-one on
$(0,\infty) \times S^1 $ and sends $\{0\} \times S^1$ to $[0]$.  The
only difference is that now $\sigma$ is given by
$$
\sigma (s, \lambda) = [\lambda^{-1/2} \sqrt{s}] .
$$

Consider now the ring of (complex valued) ``polynomial'' functions on
$\CC/\ZZ_2$.  This ring, by definition, is the ring of
$\ZZ_2$-invariant polynomial functions on $\CC$; and, as we noted
above, this ring is generated by $z^2$, $\bar{z}^2$ and $|z|^2$.  By
abuse of notation we can think of these functions as living on $\CC
/\ZZ_2$.  Now note that now $\sigma^* z^2 = \lambda\inv s =
e^{-i\theta}s$, $\sigma^* \bar{z}^2 = e^{i\theta}s$ and
$\sigma^*|z|^2=s$.  On the other hand we will prove in
\S~\ref{sec:2} that the ring of differential operators on $S^1$
which commute with $\Pi$ is generated by
\begin{displaymath}
  \frac{1}{i} \frac{d}{d\theta} e^{i\theta}, \quad
  \frac{1}{i} e^{-i\theta} \frac{d}{d\theta} \hbox{ and }
  \frac{1}{i} \frac{d}{d\theta} \, ;
\end{displaymath}
and the symbols of these operators are exactly the functions
$se^{i\theta}$, $se^{-i\theta}$ and $s$ above.  Thus the ring of
polynomial functions on $(T^*S^1)_{++}$ is exactly the ring of
symbols of differential operators which commute with $\Pi$.\\

\section{The Szeg\"o projector on $S^1$}
\label{sec:2}

The classical Szeg\"o projector
$$
\Pi : L^2 (S^1) \to L^2 (S^1)
$$
is the orthogonal projection of the space $L^2 (S^1)$ onto the
space
$$
{\span}\, \{ e^{i n\theta} \mid \, n\geq 0\}.
$$
Our goal in this
section is to determine all differential operators on $S^1$ which
commute with $\Pi$. It is easy to check that the operators
\begin{equation} \label{eq:2.1}
\frac{1}{i}\frac{d}{d\theta},\,
\frac{1}{i}\frac{d}{d\theta}e ^{i \theta}\,\quad
e ^{-i \theta}\frac{1}{i}\frac{d}{d\theta}
\end{equation}
have this property, and we will prove that the only differential
operators that commute with $\Pi$ are sums and products of these
operators.

\begin{theorem} \label{theorem:2.1}
The algebra of differential operators on the circle $S^1$ which
commute with the Szeg\"o projector $\Pi$ is generated by the operators
(\ref{eq:2.1}).
\end{theorem}

\begin{proof}
We will first prove that
\begin{equation}\label{eq:2.2}
\left( \frac{1}{i}\frac{d}{d\theta}e ^{i \theta} \right)^k =
e^{i k \theta } p_k (\frac{1}{i}\frac{d}{d\theta})
\end{equation}
where
\begin{equation}\label{eq:2.3}
 p_k (x) = (x+1)\cdots (x+k).
\end{equation}
Assume by induction that this holds for $k-1$.  Then
\begin{equation*}
\begin{split}
\left( \frac{1}{i}\frac{d}{d\theta} e ^{i \theta} \right)^k &
=\frac{1}{i}\frac{d}{d\theta}e ^{i \theta} \left( e^{i (k-1) \theta} p_{k-1}
(\frac{1}{i}\frac{d}{d\theta}) \right) \\
&= \frac{1}{i}  e^{i k \theta }
p_{k-1} (\frac{1}{i}\frac{d}{d\theta}) \\
&= \frac{1}{i}  e^{i k \theta}
\left( \frac{1}{i}\frac{d}{d\theta} + k\right) p_{k-1}
(\frac{1}{i}\frac{d}{d\theta}) \\
&=  \frac{1}{i}  e^{i k \theta } p_{k} (\frac{1}{i}\frac{d}{d\theta}) \\
\end{split}
\end{equation*}
\hfill Q.E.D.

Now let $Q$ be a differential operator of degree $d$ which commutes
with $\Pi$ and transforms under the action $\tau$ of $S^1$ by
\begin{equation} \label{eq:2.4}
\tau_\theta^* Q = e^{i k\theta} Q \tau_\theta^*, \qquad k\geq 0.
\end{equation}
Such an operator has to be of the form $e^{ik\theta } q
(\frac{1}{i}\frac{d}{d\theta})$ for some $d$-th degree
polynomial $q(x)$.  The commutator condition $[Q, \Pi]= 0$ implies that
$$
Q e^{i m\theta} = \Pi Q e^{i m\theta}
Q \Pi e^{i m\theta} = 0
$$
for $m = -k, -k +1, \ldots, -1$, so the integers $m= -k+s$, $s = 0,
\ldots, k-1$ are roots of $q$.  Thus $p_k (x)$ divides $q(x)$; and letting
$r(x) = q(x)/p_k (x)$, one has:
\begin{equation} \label{eq:2.5}
Q = r(\frac{1}{i}\frac{d}{d\theta})\left(
(\frac{1}{i}\frac{d}{d\theta}) e^{i \theta}\right)^k
\end{equation}
by (\ref{eq:2.2}).

If $Q$ transforms under the action $\tau$ of $S^1$ by
\begin{equation} \label{eq:2.6}
\tau_\theta^* Q = e^{-i k\theta} Q \tau_\theta^*, \qquad k\geq 0,
\end{equation}
the transpose of $Q$ transforms by (\ref{eq:2.4}).  Therefore the
transpose of $Q$ has to be of the form (\ref{eq:2.5}), and $Q$ itself
of the form
\begin{equation} \label{eq:2.7}
Q =\left(
\frac{1}{i} e^{-i \theta}\right)^k 
r(\frac{1}{i} \frac{d}{d\theta}).
\end{equation}

Finally let $Q$ be any differential operator on the circle commuting
with $\Pi$.  Explicitly let
$$
Q = \sum _{r= 0} ^d f_r (\theta)
\left( \frac{1}{i} \frac{d}{d\theta}\right)^r,
$$
and let $c_{k, r} $ be the $k$th Fourier coefficient of $f_r (\theta)$. Then
\begin{equation}\label{eq:2.8}
Q = \sum _k Q_k
\end{equation}
with
$$
Q_k = e^{i k\theta} \sum _{r=0}^d c_{k,r}
\left( \frac{1}{i} \frac{d}{d\theta}\right)^r .
$$
Each of the $Q_k$'s commute with $\Pi$ and transform under the
action of $S^1$ by (\ref{eq:2.4}) or by (\ref{eq:2.6}); hence it has
to be of the form (\ref{eq:2.5}) or (\ref{eq:2.7}).  In particular
$Q_k = 0 $ for $|k| > d$; so the sum (\ref{eq:2.8}) is finite, and
every summand is in the algebra generated by the operators
(\ref{eq:2.1}).
\end{proof}

We will need in \S~3 an ``even'' variant of Theorem~\ref{theorem:2.1}
(whose proof we will omit since it is essentially the same as the
proof above).

\begin{theorem} \label{theorem:2.2}
Let $\Pi^{\even}$ be the orthogonal projection from $L^2 (S^1)$ onto
the space
\begin{equation}\label{eq:2.9}
\span \{ e^{2in\theta} \mid \, n\geq 0\} .
\end{equation}
The algebra of differential operators on $S^1$ which commute with
$\Pi^{\even}$ is generated by
\begin{equation}\label{eq:2.10}
\frac{1}{i} \frac{d}{d\theta} e^{2i\theta}, \quad
\frac{1}{i} \frac{d}{d\theta} e^{-2i\theta}, \quad
\frac{1}{i} \frac{d}{d\theta} .
\end{equation}
\end{theorem}

The symbols of these operators are a Poisson subalgebra of the algebra
of $C^\infty$ functions on $T^*S^1$, and as we saw in the introduction this
algebra can be identified with the algebra of ``polynomials'' on the
space $\CC/\ZZ_2 =  (T^*S^1)_{++} $.  This proves:
\begin{theorem}\label{theorem:2.4}
The algebra of differential operators on $S^1$ which commute with the
even Szeg\"o projector $\Pi^{\even}$ has for its symbol algebra the
algebra of polynomials on the cut space $\CC =(T^*S^1)_{+} $.
\end{theorem}

What about the algebra of differential operators which commute with
the usual Szeg\"o projector?  The same argument gives:

\begin{theorem}\label{theorem:2.5}
The algebra of differential operators on $S^1$ which commute with
the Szeg\"o projector $\Pi$ has for its symbol algebra the algebra of
polynomials on the cut space $\CC/\ZZ_2 = (T^*S^1)_{++} $.
\end{theorem}

Finally we characterize smooth functions on the cut space $(T^*S^1 )_+ = \CC$
which can be extended to smooth functions on $T^* S^1$.

\begin{theorem} \label{theorem:2.3}
A function $f\in C^\infty ((T^*S^1)_+)$ has the property that
$\sigma^* f \in C^0 ([0, \infty) \times S^1)$ is the restriction of a
smooth function on $T^* S^1$ iff the infinite jet of $f$ at 0 is even,
i.e., is invariant under $z\mapsto -z$.  Here, as before, $\sigma
(s,\lambda) = \lambda^{-1} \sqrt {s}$.
\end{theorem}

\begin{proof}
If $f\in C^\infty (\CC)$ vanishes at zero to infinite order, then
$\sigma^* f$ can be extended by zero to a smooth function on $\RR
\times S^1 = T^* S^1$.  Therefore the condition on $\sigma^* f$ to extend
is the condition on the infinite jet  of $f$ at 0.  We can write  the
the jet $j^\infty f (0)$ as
$$
        j^\infty f (0)= \sum_{n=0}^\infty \sum _{k+l = n} a_{kl} z^k \bar{z}^l
$$
for some $a_{kl} \in \C$.   Since $\sigma^* z = \lambda \inv s^{1/2}$,
$$
\sigma^*( \sum _{k+l = n} a_{kl} z^k \bar{z}^l) =
        (\sum _{k+l = n} a_{kl} \lambda^{l-k}) s^{n/2}.
$$
Since $\sigma^*f$ extends to a smooth function on $T^*S^1$ iff
$\sigma^* (     j^\infty f (0))$
has no fractional powers of $s$,  we must have
$$
j^\infty f (0)= \sum_{m=0}^\infty \sum _{k+l = 2m} a_{kl} z^k \bar{z}^l,
$$
i.e., $j^\infty f (0)$ is a power series in $z^2$, $\bar{z}^2$ and
$|z|^2$. The latter is true iff $j^\infty f (0)(z, \bar{z}) = j^\infty
f (0)(-z, -\bar{z})$.
\end{proof}

\section{The Szeg\"o projector on $\RR^n \times S^1$}
\label{sec:3}

Let $\Pi_1$ be the Szeg\"o  projector on $L^2(S^1)$ (the operator
we called $\Pi$ in \S~\ref{sec:2}).  From $\Pi_1$ one gets a
projection operator,
\begin{displaymath}
  I_{\RR^n} \otimes \Pi_1
\end{displaymath}
on $L^2 (\RR^*) \otimes L^2 (S^1)$ which extends by continuity to
a projection operator
\begin{displaymath}
  \Pi : L^2 (\RR^n \times S^1) \to L^3(\RR^n \times S^1) \, .
\end{displaymath}
Our goal in this section will be to determine the commutator of
$\Pi$ in the algebra of pseudodifferential operators on $\RR^n
\times  S^1$.  For simplicity we will only consider
pseudodifferential operators of the form
\begin{equation}
  \label{eq:3.1}
  Qf = \sum_m e^{im\theta} \int q(x,\xi,\theta ,m)
     e^{ix \cdot \xi} \hat{f} (\xi ,m) \, d\xi
\end{equation}
$\hat{f}$ being the Fourier transform of $f$:
\begin{equation}
  \label{eq:3.2}
  \hat{f} (\xi ,m) = \left( \frac{1}{2\pi}\right)^{n+1}
  \int e^{-ix\cdot \xi}  e^{-im\theta} f(x,\theta) \, dx \, d \theta
\end{equation}
and $q(x,\xi ,\theta s)$ being a classical polyhomogeneous
symbol of compact support in $x$.  We can decompose $q$ into its
Fourier modes
\begin{eqnarray}
  \label{eq:3.3}
  q(x,\xi ,\theta ,s)
      &=& \sum e^{ik\theta}  q_k (x,\xi ,s) \\
\noalign{\hbox{with }} \nonumber\\
  \label{eq:3.4}
  q_k (x,\xi ,s)
     &=& \frac{1}{2\pi}\int q(x,\xi,\theta,s) e^{-ik\theta}
         \, d \theta \, ;
\end{eqnarray}
and from (\ref{eq:3.3}) we get a corresponding decomposition of
$Q$:
\begin{equation}
  \label{eq:3.5}
  Q=\sum Q_k
\end{equation}
$Q_k$ being the operator with symbol
\begin{equation}
  \label{eq:3.6}
  q_k (x,\xi ,s) e^{ik\theta} \, .
\end{equation}
Letting $p_k (x,y,m)$ be the conormal distribution
\begin{equation}
  \label{eq:3.7}
  p_k (x,y,m) = \left( \frac{1}{2\pi} \right)^{n+1}
     \int q_k (x,\xi ,m) e^{i(x-y)\cdot \xi} \, d \xi
\end{equation}
we can, by (\ref{eq:3.1})--(\ref{eq:3.2}), write the Schwartz
kernel of $Q_k$ as a sum:
\begin{equation}
  \label{eq:3.8}
  \sum_m e^{i(k+m) \theta} e^{-im \psi} p_k (x,y,m) \, .
\end{equation}
>From (\ref{eq:3.8}) we will deduce:

\begin{lemma}
\label{lem:3.1}
For $k$ positive the Schwartz kernel of $\left[ \Pi , Q_k
\right]$ is
\begin{equation}
  \label{eq:3.9}
  \sum_{-k \leq m<0} e^{i(k+m)\theta} e^{-im\psi}
     p_k (x,y,m) \, .
\end{equation}

\end{lemma}

\begin{proof}
  The Schwartz kernel of $\Pi Q_k - Q_F \Pi$ is
  \begin{eqnarray*}
    \sum_{m+k \geq 0} e^{i(k+m)\theta} e^{-im\psi}
        p_k (x,y,m)
    - \sum_{m \leq 0} e^{i(k+m)\theta} e^{-im\psi}
       p_k (x,ym)
  \end{eqnarray*}
and this difference is the same as the finite sum
(\ref{eq:3.9}).  Similarly for $k$ negative one has
\end{proof}

\begin{lemma}
\label{lem:3.2}
The Schwartz kernel of $\left[ \Pi , Q_k \right]$ is
\begin{equation}
  \label{eq:3.10}
  \sum_{k \leq m<0} e^{im\theta} e^{-i(m-k)\psi}
     p_k (x,y,m) \, .
\end{equation}

\end{lemma}

{F}rom these results we can easily read off necessary and
sufficient conditions for $Q$ and $\Pi$ to commute.

\begin{theorem}
  \label{th:3.3}
  $Q$ and $\Pi$ commute if and only if, for all $k$,
  \begin{equation}
    \label{eq:3.11}
    q_k (x,\xi ,m) = 0
  \end{equation}
for $-|k| \leq m <0 \, .$
\end{theorem}

In particular for $k>0$ this implies that there exists a
classical polyhomogeneous symbol $q^\#_k (x,\xi ,s) $ with
\begin{equation}
  \label{eq:3.12}
  q_k (x,\xi , s) = q^\#_k (x,\xi ,s) \Pi^k_{m=1} (s+m) \, .
\end{equation}

Let $Q^\#_k$ be the pseudodifferential operator with $q^\#_k$ as
symbol.  Since $q^\#_k$ doesn't depend on $\theta$ this operator
commutes with the action of $S^1$ on $\RR^n \times S^1$ and by
(\ref{eq:3.6}) and (\ref{eq:2.2})
  \begin{equation}
    \label{eq:3.13}
    Q_k = Q^{\#}_k \left( \frac{1}{\sqrt{-1}} \,
         \frac{d}{d \theta} \right)^k \, .
  \end{equation}
Similarly
\begin{equation}
  \label{eq:3.14}
  Q_{-k} = Q^\#_{-k} \left( e^{-i\theta}
      \frac{1}{\sqrt{-1}} \, \frac{d}{d\theta} \right)^k
\end{equation}
so we have proved

\begin{theorem}
  \label{th:3.4}
A necessary and sufficient condition for $Q$ to commute with
$\Pi$ is that, for every $k$, $Q_k$ have a factorization of the
form, (\ref{eq:3.13})--(\ref{eq:3.14}), the operator $Q^\#_k$
being a classical polyhomogeneous pseudodifferential operator on
$\RR^n \times S^1$ which is $S^1$ invariant.

\end{theorem}

As another application of Lemmas \ref{lem:3.1}--\ref{lem:3.2} we
will prove:

\begin{theorem}
  \label{th:3.5}
  If the symbol $q (x,\xi,\theta,s)$ of $Q$ vanishes to
  infinite order on the set $\xi \neq 0$, $s=0$ the operator, $
 [\Pi ,Q]$ is a smoothing operator.
\end{theorem}

\begin{remark*}
  If $[\Pi ,Q]$ is a smoothing operator the operator
  \begin{displaymath}
    \Pi Q \Pi + (I-\Pi) Q (I-\Pi)
  \end{displaymath}
differs from $Q$ by a smoothing operator and commutes with
$\Pi$.  In other words $Q$ is the sum of an operator which
commutes with $\Pi$ and a smoothing operator.
\end{remark*}

\begin{proof}
  It suffices to show that each of the operators $[\Pi , Q_k]$ is
  smoothing and hence, by (\ref{eq:3.9}), that $p_k (x,y,m)$ is
  smooth.  But $p_k (x,y,m)$ is defined by the integral
  (\ref{eq:3.7}), and we can expand the integrand in a finite
  Taylor series
  \begin{eqnarray*}
    q_k (x,\xi ,m)
       &=& \sum^N_{\ell =0} \frac{1}{\ell !}
              \left( \frac{d}{ds} \right)^{\ell} q_k
              (x\xi ,0) m^{\ell}
        + r_N (x,\xi ,m)\\
\noalign{\hbox{where}}\\
    r_N (x,\xi ,s)
       &=& \frac{1}{N!} \int^1_0 (1-t)^N
             \left( \frac{d}{ds} \right)^N q_k
               (x,\xi ,ts) \, dt
  \end{eqnarray*}
is a classical polyhomogeneous symbol of degree equal to $\deg
Q-N$.  Thus if $q_k (x;\xi ,s)$ vanishes to infinite order at
$s=0$
\begin{displaymath}
  q_k (x,\xi ,m) = r_N (x,\xi ,m)
\end{displaymath}
for all $N$; so by (\ref{eq:3.7})
\begin{displaymath}
  p_k (x,y,m) = \left( \frac{1}{2n}\right)^{n+1}
     \int r_N (x,\xi ,m) e^{i(x-y) \cdot \xi} \, d \xi
\end{displaymath}
and, for all integers, $\ell$, the right side is in $C^{\ell}$
for $N \geq n + \deg Q + \ell$.  Hence the left hand side is in $\Cinf$.

\end{proof}
Let $Q$ be a pseudodifferential operator of order $m$ which
commutes with $\Pi$, and let $\sigma = \sigma (Q) (x,\xi,\theta
,s)$ be its leading symbol.  By (\ref{eq:3.4}) this leading
symbol only depends on the variables $s$ and $\theta$, as a
smooth function of $s$, $se^{i\theta}$ and $se^{-i\theta}$.  We
will prove that the converse is true.

\begin{theorem}
  \label{th:3.6}
  Let $\sigma$ be a smooth function on the complement of the zero
  section in $T^* (\RR^n \times S^1)$ which is homogeneous of
  degree $m$ and only depends on $s$ and $\theta$ as a smooth
  function of $s$, $se^{i\theta}$ and $se^{-i\theta}$.  Then
  there exists an $m$\st{th} order pseudodifferential operator, $
 Q$, which commutes with $\Pi$ and has leading symbol, $\sigma$.
\end{theorem}

\begin{proof}
  Let $\sigma = \sigma_+ + \sigma_- +\sigma_0$, $\sigma_+$ being
  the sum of the positive Fourier modes of $\sigma$ and
  $\sigma_-$ the sum of the negative Fourier modes.  Since
  $\sigma_-$ is the complex conjugate of $\bar{\sigma}_+$, it
  suffices to prove the theorem for $\sigma_+$.  Let $\sigma_k$,
  $k>0$, be the $k$\st{th} Fourier mode of $\sigma$.  By
  hypothesis
  \begin{displaymath}
    \sigma_k (x,\xi ,\theta ,s)
       = \sigma^{\#}_k (x,\xi,s)s^k e^{ik\theta} \, .
  \end{displaymath}
Let $Q^\#_k$ be an $S^1$ invariant pseudodifferential operator
with leading symbol equal to $\sigma_k$, and let
\begin{displaymath}
  Q_k = Q^\#_k \left( \frac{1}{\sqrt{-1}} \,
    \frac{d}{d\theta} e^{i\theta} \right)^k \, .
\end{displaymath}
Then $Q_k$ commutes with $\Pi$ and has $\sigma_k$ as its leading
symbol.  Let $H$ be the pseudodifferential operator on $\RR^n
\times S^1$ with symbol $(\xi^2 + s^2)^{-\frac{1}{2}}$, let $\rho
(s)$ be a compactly supported function which is $1$ on the
interval, $|s| <1$ and let
\begin{displaymath}
  N_1 < N_2 < \cdots
\end{displaymath}
be an increasing sequence of positive integers.  Then the sum
\begin{displaymath}
  Q_+ = \sum_{k>0} \rho \left( N_k \left(
      \frac{1}{\sqrt{-1}} \, \frac{d}{d\theta}
    \right) H \right) Q_k
\end{displaymath}
is well-defined and (provided that the $N_k$'s go to infinity
fast enough) is a classical pseudodifferential operator which
commutes with $\Pi$ and has leading symbol
\begin{displaymath}
  \sum_{k>0}\rho \left( N_k \frac{s}{(\xi^2 + s^2)^{\frac{1}{2}}}
  \right)
  \sigma_k (x,\xi,s,\theta) \, .
\end{displaymath}
In particular this symbol has the same formal power series
expansion on the set $s=0$ as does $\sigma_+$.  Hence one can
find an $m$\st{th}  order pseudodifferential operator, $R_+$,
whose total symbol vanishes to infinite order on $s=0$ and whose
leading symbol is $\sigma_+ -\sigma (Q_+)$.  Thus $Q_+ + R_+$ has
leading symbol, $\sigma_+$, and commutes with $\Pi$ modulo
smoothing operators.  Therefore, as we pointed out above, it is
the sum of an operator which commutes with $\Pi$ and a smoothing
operator.

\end{proof}

Let $Q$ be a pseudodifferential operator which commutes with
$\Pi$.  We will show that the operator
\begin{displaymath}
  \Pi Q = Q \Pi = \Pi Q \Pi
\end{displaymath}
``lives microlocally'' on the set $s>0$.

\begin{theorem}
  \label{th:3.7}
  $\Pi Q$ is smoothing if and only if the symbol $q(x,\xi,\theta
  ,s)$ of $Q$ is of order $-\infty$ on the set $s \geq 0$.
\end{theorem}

\begin{proof}
  By (\ref{eq:3.1})
  \begin{displaymath}
    \Pi Q f = \sum_{m \geq 0} e^{im\theta} \int q
    (x,\xi,\theta,m) e^{ix \cdot \xi} \hat{f} (\xi , m) \, d\xi
  \end{displaymath}
and this is smoothing if and only if $q$ is a symbol of order
$-\infty$ on the set $s \geq 0$.

Let $(\Pi_1)_{\even}$ be the even Szeg\"o projector on
$L^2(S^1)$ (the operator we called $\Pi_{\even}$ in
\S~\ref{sec:2}) and let
\begin{displaymath}
  \Pi_{\even} = I_{\RR^n} \otimes (\Pi_1)_{\even} \, .
\end{displaymath}
For this projector there are obvious analogues of
Theorems~(\ref{th:3.3})--(\ref{th:3.7}).  We will content ourselves with
describing the even analogue of Theorem~\ref{th:3.4}.

\end{proof}

\begin{theorem}
  \label{th:3.8}
 A necessary and sufficient condition for $Q$ to commute with
$ \Pi_{\even}$ is that, for all $k$, $Q_{2k+1} =0$, and for all
 positive $k$
 \begin{eqnarray}
\label{eq:3.15}
   Q_{2k} &=& Q^\#_{2k} \left( \frac{1}{\sqrt{-1}} \,
     \frac{d}{d\theta} e^{2i\theta}\right)^k \\
\noalign{\hbox{and}}\\
\label{eq:3.16}
Q_{-2k} &=& Q^\#_{-2k} \left( e^{-2i\theta} \frac{1}{\sqrt{-1}} \,
  \frac{d}{d\theta}\right)^k
 \end{eqnarray}
$Q^\#_{2k}$ and $Q^\#_{-2k}$ being pseudodifferential operators
which are $S^1$-invariant.
\end{theorem}

\section{Canonical forms for circle actions}\label{sec:4}

The first of the canonical forms which we will discuss in this section
is an equivariant Darboux theorem for symplectic cones. We recall that
a symplectic cone is a symplectic manifold $(M, \omega)$ equipped with
a free proper action $\rho$ of $\RR$ which satisfies
\begin{equation}\label{eq:4.1}
{\rho_a}^* \omega = e^a\omega
\end{equation}
Let $\Xi$ be the vector field generating the action, $\Xi (m) =
\left. \frac{d}{dt}\right|_{t=0} \rho_t (m)$.  The infinitesimal version of
(\ref{eq:4.1}) is
\begin{equation}\label{eq:4.2}
\omega = L_\Xi \omega = d(\iota(\Xi) \omega).
\end{equation}

Suppose now that in addition to the $\RR$ action one has a free action
$\tau$ of $S^1$ on $M$ which preserves the symplectic form $\omega$
and commutes with $\rho$, hence preserves
\begin{equation} \label{eq:4.4}
\alpha := \iota (\Xi)\omega.
\end{equation}
Then, if we denote the generator of the $S^1$ action by $V$,
$$
0= L_V \alpha = \iota (V) d\alpha + d\iota (V)\alpha .
$$
Since $\omega = d\alpha$, we get
\begin{equation}\label{eq:4.5}
\iota(V) \omega = -d (\alpha (V)).
\end{equation}
In other words $\tau$ is a Hamiltonian action with moment map
\begin{equation} \label{eq:4.6}
\Phi = \alpha (V) .
\end{equation}

Let $d = \dim M/2 = n+1$.  A simple canonical model for a $2d$
dimension symplectic cone with a homogeneous symplectic action of
$S^1$ is the complement $M_0$ of the zero section in $T^* (\RR^n \times
S^1)$.  In this model
$$
\omega_0 = \sum d\xi_i \wedge dx_i + ds \wedge d\theta
$$
is the symplectic form,
$$
\alpha _0  = \sum \xi_i dx_i + s  d \theta
$$
is the Liouville one-form (so that   $d\alpha_0 = \omega_0$),
$$
\Xi_0 = \sum \xi_i \frac{\partial}{\partial \xi_i} +
s \frac{\partial}{\partial \theta}
$$
is the generator of the $\RR$ action (so that
$\iota (\Xi_0) \omega_ 0 = \alpha _0$),
$$
V_0 = \frac{\partial}{\partial \theta}
$$
is the generator of the $S^1$ action and
$$
\Phi_0 = s
$$ is the corresponding moment map.

\begin{theorem} \label{theorem:4.1}
Let $(M, \omega, \Phi: M \to \RR)$ and $(M_0, \omega_0, \Phi_0 : M_0
\to \RR)$ be as above.  Let $p$ and $p_0$ be points in $M$ and $M_0$
respectively.  If $\Phi (p) = \Phi_0 (p_0)$, there exist $S^1 \times
\RR$ invariant neighborhoods $U$ and $U_0$ of $p$ and $p_0$
respectively, and an $S^1 \times \RR$ equivariant symplectomorphism $\gamma$ of
$(U, p)$ onto $(U_0, p_0)$.
\end{theorem}

\begin{proof}
Let $\xi = \Xi (p)$, $v= V(p)$, $\xi_0 =\Xi_0 (p_0)$ and $v_0 =
v(p_0)$. By definition of $\alpha$ (equation (\ref{eq:4.4}))
$$
\omega_p (\xi, v) = \alpha _p (v) = \Phi (p)
$$
and
$$
(\omega_0)_p (\xi_0, v_0) = (\alpha_0) _p (v_0) = \Phi_0 (p_0)
$$
so
$$
\omega_p (\xi, v) = (\omega_0)_p (\xi_0, v_0).
$$
Hence there exists a linear symplectic mapping $A:T_pM \to T_{p_0}
M_0$ mapping $\xi$ to $\xi_0$ and $v$ to $v_0$ (note that there are
two cases to consider: $\omega_p (\xi, v) = 0$ and $\omega_p (\xi, v)
\not =0$).  Let $X$ and $X_0$ be the $S^1\times \RR$ orbits through $p$
and $p_0$, and $i$ and $i_0$ the inclusions of $X$ and $X_0$ into $M$
and $M_0$ respectively.  The map $A$ above extends uniquely to an $S^1
\times \RR$ equivariant isomorphism of symplectic vector bundles
$$
A: i^* TM \to i_0^* TM_0,
$$
and this can be exponentiated to an $S^1 \times \RR$ equivariant map
$$
\Gamma_A : W \to W_0
$$
 of an $S^1 \times \RR$ neighborhood $W$ of $X$ onto an $S^1 \times
\RR$ invariant neighborhood $W_0$ of $X_0$ with the property that
$d\Gamma_A = A$ at the points of $X$.  Indeed, since the action of $S^1
\times \RR$ is proper, the orbits $X$ and $X_0$ are embedded, and there
is an $S^1 \times \RR$-invariant metrics on $M$ and $M_0$.  Use the
exponential maps for these invariant metrics.  By construction of
$\Gamma_A$, the form $\tilde{\omega}: = \Gamma_A^* \omega_0$ is equal
to $\omega$ at all points of $X$.  To conclude the proof of the
theorem, we will show that there exists an $S^1 \times \RR$ invariant
neighborhood $U$ of $X$ and an $S^1 \times \RR$ equivariant open
embedding $f: U \to W$ such that $f = id$ on $X$ and
$f^*\tilde{\omega} = \omega$.  The proof will be the standard Moser
deformation argument. However we must check that it produces an $S^1
\times \RR$ equivariant deformation.

Let $\omega_t =(1-t) \omega + t\tilde{\omega}$ and let $W'$ be the open
subset of $W$ on which $\omega_t$ is symplectic for $0\leq t \leq 1$.
Since $\omega_t$ is $S^1 \times \RR$ invariant, $W'$ is $S^1 \times \RR$
invariant.  Since $\omega_t = \omega$ on $X$, $W'$ contains $X$.
Since $\omega_t$ is non-degenerate on $W'$ there exists a vector field
$y_t$ on $W'$ satisfying
\begin{equation} \label{eq:4.7}
\iota(y_t) \omega_t = \alpha - \tilde{\alpha},
\end{equation}
where $\tilde{\alpha} := \iota (\Xi) \tilde{\omega}$.  Moreover, since
$\alpha - \tilde{\alpha} =\iota(\Xi) (\omega - \tilde{\omega})$, the
vector field $y_t$ is zero at points of $X$.  Let $U $ be the subset
of $W'$ consisting of all the points $q$ at which $y_t$ has an
integral curve $\gamma_q (t)$ with $\gamma_q (0) =q$ and $\gamma_q
(t)$ is defined for all $t\in [0,1]$. Let $ f_t : U\to W' $ be the map
$f_t (q) = \gamma_q (t)$, the isotopy generated by $y_t$. Then by
Moser's trick $f_t^* \omega_t = \omega$, and in particular $f_1^*
\tilde{\omega} = \omega$.  Note that by definition $U$ and $W'$ are
$S^1 \times \RR$ invariant.  Moreover, by (\ref{eq:4.1}) $\rho_a^*
\omega_t = e^a \omega_t$ and $\rho_a^*(\alpha - \tilde{\alpha}) = e^a
(\alpha - \tilde{\alpha})$; so the vector field $y_t$ defined by
(\ref{eq:4.7}) is $\RR$ invariant. Thus the isotopy it generates is
$\RR$ equivariant.
\end{proof}

The second canonical form is for $S^1$ representations quantizing
canonical actions.  Let $M$ be the cotangent bundle of a compact manifold $X$ with the
zero section removed: $M = T^*X \smallsetminus X$, and 
$\tau$ be an action of $S^1$ on $M$ which preserves the
canonical  one form.  Let
$X_0 = \RR^n \times S^1$ ($n = \dim X -1$), $M_0 = T^*X_0
\smallsetminus X_0$, and $\tau_0$ the obvious action of $S^1$ on
$M_0$.

These actions quantize to give representation $\tau^{\#}$ and $\tau^{\#}_0$, of $S^1$ on $L^2
(X)$ and $L^2 (\RR^n \times S^1)$. 
Let $p$ and $p_0$ be points of $M$ and $M_0$ with $\Phi (p) =
\Phi_0 (p_0)$, where $\Phi$ and $\Phi_0$ are the corresponding moment
maps, and  let $\gamma : (U, p) \to (U_0, p_0)$ be a canonical
transformation mapping an $S^1 \times \RR$ -invariant neighborhood $U$
of $p$ onto an $S^1 \times \RR$ -invariant neighborhood $U_0$ of $p_0$.

\begin{theorem} \label{theorem:4.2}
The transformation $\gamma$ can be implemented by a Fourier integral
operator of order zero
$$F : \Cinf (X) \to \Cinf (\RR^n \times S^1)$$
with the properties
\begin{enumerate}
\item $F^* F = I$ on $U$,
\item $FF^* = I $ on $U_0$,
\item $\tau^\# _0 (e^{i\theta}) F = F \tau^\#  (e^{i\theta})$.
\end{enumerate}
\end{theorem}
\begin{proof}
Let $F_0$ be the zeroth order Fourier integral operator with compact
support which implements $\gamma$ on $U$ and has the following three
properties.
\begin{enumerate}
\item $F_0 ^* F_0 = I + R_0 $ on $U$, $R_0$ being a pseudodifferential
operator of order $-1$.
\item $F_0 F_0 ^* = I + S_0$ on $U_0$, $S_0$ being a pseudodifferential
operator of order $-1$.
\item The symbol of $F_0$ is $S^1$ invariant.
\end{enumerate}
By averaging $F_0$ by the action
$$
\theta \mapsto \tau_0^\# (e^{-i\theta}) F_0  \tau_0^\# (e^{i\theta})
$$ one gets a Fourier integral operator $F_1$ which implements
$\gamma$, has the same leading symbol as $F_0$ and intertwines
$\tau^\#$ and $\tau_0^\#$.  In particular, since it has the same a
leading symbol as $F_0$ it continues to satisfy $F_1 ^* F_1 = I + R $,
$F F ^* = I + S$ with pseudodifferential operators of order $-1$ $R$
and $S$. Now define $F$ to be the operator
$$
F_1 (I+R)^{-\frac{1}{2}} = (I+S)^{-\frac{1}{2}} F_1 .
$$
\end{proof}

Let $\Pi$ and $\Pi_0$ be the Szeg\"o projections associated with the
representations $\tau^\#$ and $\tau_0^\#$.  One consequence of
Theorem~\ref{theorem:4.2} is that if $Q$ is a pseudodifferential
operator with microsupport on $U$ which commutes with $\Pi$ modulo
smoothing operators, $FQF^*$ is a pseudodifferential operator with
microsupport in $U_0$ which commutes with $\Pi_0$ modulo smoothing
operators.  Hence many of the results which we proved in \S~3 for the
commutator ring of $\Pi_0$ are valid for the commutator ring of $\Pi$
as well.  We will describe a number of such results in the next
section.

\section{The algebra, $\Psi^+$, and its symbol calculus}
\label{sec:5}

Let $M$ be the cotangent bundle of a compact manifold, $X$, with
its zero section deleted, let $\tau$ be an action of $S^1$ on $M$
by canonical transformations, and let $\tau^\#$ be a
representation of $S^1$ on $L^2 (X)$ compatible with $\tau$.  Let
$\Pi^{\even}$ be the ``even'' Szeg\"o projector (defined in
\S~\ref{sec:intro}) and $\Psi^{\even}_+$ the algebra of
$\ZZ_2$-invariant pseudodifferential operators on $M$ which commute with
$\Pi^{\even}$.  As we pointed put in the introduction, the
complement, $U$, of $M_{\red}$ in $M_+$ can be identified with the
open set, $U$ in $M$ where the moment map of $\tau$ is positive.
Therefore, if $A$ is a pseudodifferential operator of order $m$
in $\Psi^{\even}_+$, the restriction of its leading symbol to $U$
can be regarded as a homogeneous function of degree $m$ on the
open dense subset, $U$, of $M_+$.  By Theorems~\ref{theorem:2.4},
\ref{theorem:4.2} and the even version of Theorem~\ref{th:3.6}, this
function extends to a smooth even function on $M_+$.  Thus one
has a symbol map
\begin{equation}
  \label{eq:5.1}
  \left( \Psi^{\even}_+ \right)_m \to \Cinf (M_+)^{\even}_m
\end{equation}
from the space of $m$\st{th} order pseudodifferential operators in
$\Psi^{\even}_+$ to the space of even homogeneous functions of
degree $m$ on $M_+$.  Let $\A$ be the algebra of operators
\begin{displaymath}
  \Pi^{\even} \Psi^{\even}_+ \Pi^{\even}
\end{displaymath}
and let
\begin{displaymath}
  \A^m = \Pi^{\even} \left( \Psi^{\even}_+ \right)_m \Pi^{\even}
  \, .
\end{displaymath}

\begin{theorem}
  \label{th:5.1}
From the map (\ref{eq:5.1}) one gets a short exact sequence
\begin{equation}
  \label{eq:5.2}
  0 \to \A^{m-1} \to \A^m \overset{\sigma}{\longrightarrow}
  \Cinf (M_+)^{\even}_m \to 0 \, .
\end{equation}

\end{theorem}

\begin{proof}
  Given $A \in (\Psi^{\even}_+)_m$, suppose the leading symbol of
  $A$ vanishes on $U$.  Then one can find a pseudodifferential
  operator, $A'$, whose total symbol vanishes on $U$ and whose
  leading symbol is identical with $\sigma (A)$.  Thus by
  Theorem~\ref{th:3.5} $A'$ commutes with $\Pi^{\even}$ modulo
  smoothing operators.  Hence, by the remark following
  Theorem~\ref{th:3.5}, one can modify $A'$ by adding to it a
  smoothing operator, so it actually does commute with
  $\Pi^{\even}$.  Moreover, since the total symbol of $A'$
  vanishes on $U$ $\Pi^{\even} A'$ is smoothing by the even
  version of Theorem~\ref{th:3.7}
; so by replacing $A'$ by $A' - \Pi^{\even A'}$, one can assume
not only that $A'$ commutes with $\Pi^{\even}$ but that
$\Pi^{\even} A' =0$.  Since $\sigma (A) = \sigma (A')$, the operator,
$A-A'$ is of order $m-1$ and
\begin{displaymath}
  \Pi^{\even} A \Pi^{even} = \Pi^{\even} (A-A') \Pi^{\even} \, .
\end{displaymath}
This proves that the map, $\sigma$, in (\ref{eq:5.2}) is
injective; and that it is surjective follows from (the even
version of) Theorem~\ref{th:3.6}.

\end{proof}

We  claim next

\begin{theorem}
  \label{th:5.2}
  If $A_1$ and $A_2$ are in $\A$, $\sigma (A_1A_2) = \sigma (A_1)
  \sigma (A_2)$ and $\sigma ([A_1,A_2]) = - \sqrt{-1} \{ \sigma
  (A_1),\sigma (A_2) \}$.  Moreover if $A \in \A$, $A^* \in \A$
  and $\sigma (A^*) = \overline{\sigma (A)}$.
\end{theorem}

\begin{proof}
  Microlocally on $U$ these are standard identities for leading
  symbols of pseudodifferential operators.  Therefore, since $U$
  is a dense subset of $M_+$ they hold globally on all of $M_+$.
\end{proof}

An operator, $A \in \A^m$ is \emph{elliptic} if $\sigma (A)$ is
everywhere non-zero.  We will show that these operators have the
usual properties of elliptic operators:

\begin{theorem}
  \label{th:5.3}
  If $A \in \A^m$ is elliptic, it is invertible modulo smoothing
  operators, i.e.,~there exists a $B \in \A^{-m}$ such that $I-BA$
  and $I-AB$ are smoothing.
\end{theorem}

\begin{proof}
  Replacing $A$ by $A^*A$ we can assume that $A$ is self-adjoint
  and that $\sigma (A) >0$.  Let $A =\Pi^{\even} Q \Pi^{\even}$,
  $Q \in \Psi^{\even}_+$.  Since $\sigma (Q) = \sigma (A)$ on $U$
  we an assume that $\sigma (Q)>0$ on an open conic set, $V$ in $
 M$ containing the closure of $U$.  Let $P$ be a
 pseudodifferential operator of order $m$ whose total symbol is
 supported in the complement of the closure of $U$ and whose
 leading symbol is non-negative and strictly greater than zero on
 the complement of $V$.  By Theorems~\ref{th:3.5} and \ref{th:3.6}
   $[P,\Pi^{\even}]$ and $\Pi^{\even}P$ are smoothing, so, by
   modifying $P$ by a smoothing operator, we can assume that
   $[P,\Pi^{\even}]$ and $\Pi^{\even}P$ are zero.  Replacing $Q$
   by $Q+\lambda P$, $\lambda \gg 0$, we can assume that the
   symbol of $Q$ is positive everywhere, and hence that $Q$ is
   invertible modulo smoothing operators, i.e.,~there exists a
   pseudodifferential operator, $Q_1$, of order $-m$, with
   $Q_1Q-I$ and $QQ_1 -I$ smoothing.  It is easy to see that
   $[\Pi ,Q_1]$ is smoothing; and hence $Q_1$ can be modified by
   adding to it a smoothing operator such that $[\Pi ,Q_1]=0$.
   Now set $B= \Pi Q_1 \Pi$.
\end{proof}

The results above justify to some extent the assertion in
Theorem~\ref{th:5.1} that the algebra $\A$ ``quantizes'' the
algebra of classical observables, $\Cinf (M_+)^{\even}$.  A
slightly more compelling justification is the following.

\begin{theorem}
  \label{th:5.4}
  If $A \in \A_m$, $m >0$, is elliptic and self-adjoint and
  $\sigma (A)$ is everywhere--positive, the spectrum of $A$ is
  discrete, and its eigenvalues
  \begin{displaymath}
    \lambda_1 \leq \lambda_2 \leq \cdots
  \end{displaymath}
satisfy the Weyl law
\begin{displaymath}
  N(\lambda) \sim \vol \{ m \in M_+ , \sigma (A) (m) < \lambda \}
  \, .
\end{displaymath}
Here $N (\lambda)$ is the Weyl counting function
\begin{displaymath}
  N(\lambda) = \# \{ \lambda_i < \lambda \}
\end{displaymath}
and ``vol'' means symplectic volume.

It is shown in \cite{Gu2} that a Weyl law for an algebra of
operators of the type above is implied by the existence of a
``residue trace''; and the following theorem asserts that a
``residue trace'' exists on the algebra $\A$.

\end{theorem}

\begin{theorem}
  \label{th:5.5}
There exists a linear map
\begin{displaymath}
  \res : \A \to \CC
\end{displaymath}
with the following properties.

\alphaparenlist
\begin{enumerate}
\item 
  $\res A =0$ if and only if $A$ can be written as a sum of
  commutators
  \begin{displaymath}
    A=\sum^N_{i=1} [A_i,B_i] \, ,
  \end{displaymath}
$A_i,B_i \in \A$.

\item 
  If $A$ is of degree $-n$
  \begin{equation}
    \label{eq:5.3}
    \res (A) = \int_{M_+} \sigma (A) \omega^m_+
  \end{equation}
$\omega_+$ being the symplectic form on $M_+$.

\end{enumerate}

\end{theorem}

\begin{remark*}
  If $(M,\omega)$ is a symplectic cone of dimension $2n$, and $f$
  a homogeneous function of degree $-n$, the form $f \omega^n$ is
  a $2n$ form of degree of homogeneity zero; so
  \begin{displaymath}
    L_{\Xi} f\omega^n =0=d (\iota (\Xi) f\omega^n) \, .
  \end{displaymath}
Thus $\iota (\Xi) f\omega^n$ is closed.  Let $\Gamma$ be a
compact $2n-1$ dimensional submanifold of $M$.  $\Gamma$ is
called a \emph{contour} if it intersects every ray of the cone,
$M$, in exactly one point.  It is very easy to see that if
$\Gamma$ and $\Gamma_1$ are contours, $\Gamma$ can be smoothly
deformed into $\Gamma_1$ and hence the integral
\begin{displaymath}
  \int_{\Gamma} \iota (\Xi) f\omega^n
\end{displaymath}
is independent of the choice of $\Gamma$; and this integral is
\emph{defined} to be the integral
\begin{displaymath}
  \int f\omega^n \, .
\end{displaymath}

\end{remark*}

We won't give the proof of the existence of this residue trace
here.  Details can be found in \cite{Gu}.

We will next describe some analogous results for the Szeg\"o
projector, $\Pi$, and the algebra of the pseudodifferential
operators, $\Psi_+$, commuting with $\Pi$.  Let
\begin{displaymath}
  \B = \Pi \Psi_+ \Pi \, .
\end{displaymath}
%
%
It is clear from Theorems~\ref{theorem:2.5} and \ref{th:3.6} (and
the canonical form Theorem~\ref{theorem:4.2}) that the leading
symbol of an operator, $B \in \B^n$ can be interpreted as a
function on $M_{++}$; and the following is proved by the
same proof as that of Theorem~\ref{th:5.1}.

\begin{theorem}
  \label{th:5.6}
  There exists a short exact sequence
  \begin{equation}
    \label{eq:5.4}
    0 \to \B^{m-1} \to \B^n \overset{\sigma}{\longrightarrow}
      \Cinf (M_{++})_m \, .
  \end{equation}

\end{theorem}

Notice, by the way, that if $\tilde{U}$ is the complement of the
cut locus, $M_{\red}$, in $M_{++}$, one has a map
\begin{equation}
  \label{eq:5.5}
  \gamma : \B^m \to \Cinf (M_{++})_{\even}
        \to \Cinf (U) \, .
\end{equation}
We claim:

\begin{theorem}
  \label{th:5.7}

If $B \in \B^m$ is of the form $B=\Pi Q \Pi$, with $Q \in
\Psi^m_+$, then $\gamma (B)$ is the restriction to $U$ of the
usual pseudodifferential symbol of $Q$.
\end{theorem}

In other words on the complement of the cut locus in $M$
the symbol calculus for the algebra, $\B$, is identical with the
usual symbol calculus for pseudodifferential operators on the
open subset, $U$, of $M$.

\begin{proof}
  It suffices to check this in the model case, $X = \RR^n \times
  S^1$; and in this model case, it is a consequence of
  Theorem~\ref{th:3.5} and Theorem~\ref{th:3.7}.
\end{proof}

\section{Application: toric symplectic cones} \label{sec:6}

In this section we apply our microlocal version of symplectic cuts to
the punctured cotangent bundle $T^*_0 S^2 := T^*S^2 \smallsetminus
S^2$ of the two-sphere to obtain symplectic cones over lens spaces.
We then show that by applying symplectic cuts repeatedly to the
punctured cotangent bundle of an $n$-torus one can obtain almost all
symplectic toric cones.

As a preparation for the argument to follow, we generalize
Proposition~\ref{prop:cut}(1) (see \cite{Le} for details).  Suppose
we have a Hamiltonian action of an $n$-torus $G\simeq \RR^n/\ZZ^n$ on
a symplectic manifold $(M, \omega)$ with an associated moment moment
map $\Phi :M \to \fg^*$.  Pick a primitive vector $\lambda$ in the
integral lattice $\ZZ_G$ of $G$.  Then the group $H_\lambda:= \{ \exp
t \lambda \mid t\in \RR\}$ is a closed subgroup of $G$ isomorphic to
$S^1$.  The restriction of the action of $G$ on $M$ to $H\lambda$ is
Hamiltonian with a corresponding moment map
$$
\Phi_\lambda = \langle \Phi , \lambda \rangle ,
$$
where, as usual, $ \langle , \rangle : \fg^* \times \fg \to \RR$ is
the canonical pairing.  If $H_\lambda $ acts freely on the set
$\Phi_\lambda \inv (0)$, then the cut of $M$ with respect to the
action of $H_\lambda$ makes sense.  We denote the resulting space by
$M_{+\lambda}$.  Since the actions of $G$ and $H_\lambda$ on $M$
commute, the action of $G$ on $M$ descends to a Hamiltonian action of
$G$ on $M_{+\lambda}$. The moment map $\Phi$ descends to a map
$\Phi_{+\lambda}$ on $ M_{+\lambda}$; it is an associated moment map
for the action of $G$.  Finally, it is not hard to see that
$$
 \Phi_{+\lambda} (M_{+\lambda}) = \Phi (M) \cap \{ \eta \in\fg^* \mid
\langle \eta, \lambda \rangle \geq 0\}.
$$
In other words the moment image of $M_{+\lambda}$ is cut out from the
moment image of $M$ by the half-space $\{\eta \mid \langle \eta,
\lambda \rangle \geq 0\}$.

Another ingredient that we will need is an analogue of the Delzant's
theorem for toric symplectic cones.  Recall that a toric symplectic
cone is a symplectic manifold $ (M, \omega)$ with a free proper action
$\{\rho_t\}$ of $\RR$ making it a symplectic cone and with an
effective symplectic action of a torus $G$ commuting with $\{\rho_t\}$
and satisfying $2\dim G = \dim M$. (Note that such an action of $G$ is
automatically Hamiltonian and that there is a naturally associated
moment map $\Phi : M \to \fg^*$ with $\Phi (\rho_t (m)) = e^t \Phi
(m)$ for all $m\in M$, $t\in \RR$.)\,\, We will further assume
throughout that the base $M/\RR$ of our symplectic cone $(M, \omega,
\rho_t, \Phi :M \to \fg^*)$ is compact and connected.  Note that the
base $M/\RR$ is naturally contact; more or less by definition it is a
contact toric manifold.

\begin{remark}\label{last_rem}
The classification of compact connected contact toric manifolds
(equivalently, of symplectic toric cones over a compact connected
base) is somewhat more complicated than Delzant's classification of
compact symplectic toric manifolds; see \cite{CTM} and references
therein.  There is, however, a class of symplectic toric cones for
which the classification is particularly nice.  Namely assume
in addition, the moment image $\Phi (M)$ lies in an open
half-space in $\fg^*$, i.e., that there is a vector $X\in \fg$ such
that the function $\langle \Phi, X \rangle $ is strictly positive.
Then $\Phi (M) \cup \{0\}$ is a strictly convex rational\footnote{
``rational'' means that the supporting hyperplanes are cut out by
vectors in the integral lattice of the torus $G$} polyhedral cone
(the result is implicit in \cite{BG}; cf.\ \cite[Theorem~4.3]{LS}).
Moreover, the polyhedral cone $\Phi(M) \cup \{0\}$ uniquely determines
the symplectic toric cone $(M, \omega, \rho_t, \Phi :M \to \fg^*)$.
In particular, if $(M_i, \omega_i, \rho^i_t, \Phi_i :M_i \to \fg^*)$,
$i=1,2$, are two symplectic toric $G$-cones (over a compact connected
base) whose moment images are the same convex polyhedral cones, then
$M_1$ and $M_2$ are isomorphic as symplectic toric $G$-cones
\cite[Theorem~2.18 (4)]{CTM}.
\end{remark}

In what follows we take the standard $n$ torus $\bbT^n$ to be the Lie
group $\RR^n/\ZZ^n$.  Thus the Lie algebra of $\bbT^n$ is $\RR^n$.
The identification of $\RR^n$ with $(\RR^n)^* $ by way of the standard
basis identifies the weight lattice of $\bbT^n$ with $\ZZ^n$.

Consider the action of the torus $\bbT^2$ on the punctured cotangent
bundle $T^*_0 S^2$ generated by the normalized geodesic flow for the
round metric and by the lift of a rotation of $S^2$ about an axis.  It
is not hard to see that the image of the associated homogeneous moment
map
$$
\Phi : T^*_0 S^2 \to \RR^2
$$
is the cone $C$ spanned by the vectors $(-1, 1)$ and $(1, 1)$ with
the vertex at the origin deleted:
$$
C = \{ t_1 (-1, 1) + t_2 (1, 1) \in \RR ^2 \mid t_1, t_2 \geq 0\}
\, ,
$$
so the manifold $T^*_0 S^2$ is a symplectic cone over $\R P^3$.

More generally there is a natural action of $\bbT^2$ on the symplectic
cone over any lens space $L(p, q)$.  Fix two positive relatively prime
integers $p$ and $q$.  The map $\bbT^2 \to S^1 \times S^1 = \{ (z_1,
z_2 ) \in \C^2 \mid |z_1|^2 = |z_2|^2 = 1\}$, $[\theta_1, \theta_2]
\mapsto (e^{2\pi i\theta_1}, e^{2\pi i \theta_2})$ identifies $\bbT^2
$ with $S^1 \times S^1$.  The group $\Gamma := \{(\mu_1, \mu_2 ) \in
S^1 \times S^1 \mid \mu_1 \mu_2 ^q = 1, \, \mu_1 ^p = 1\}$ is cyclic
of order $p$.  The quotient of $\C^2 \smallsetminus \{0\}$ by the
natural action of $\Gamma$ ($(\mu_1, \mu_2) \cdot (z_1, z_2) = (\mu_1
z_1, \mu_2 z_2 )$) is, more or less by definition, the symplectic cone
on the lens space $L(p, q)$:
$$
(\C^2 \smallsetminus \{0\})/\Gamma = L(p, q) \times \RR .
$$
The natural action of $\bbT^2 \simeq S^1 \times S^1$ on $\C^2$
descends to an effective Hamiltonian action of $\bbT^2 /\Gamma$ on the
cone $ L(p, q) \times \RR$.  We compute the image of the associated
moment map as follows.  The natural action of $\bbT^2 \simeq S^1
\times S^1$ on $\C^2$ descends to an effective Hamiltonian action of
$\bbT^2/\Gamma$ on $\C^2 \smallsetminus \{0\})/\Gamma$.  The kernel of
the surjective map $\varphi : S^1 \times S^1 \to S^1 \times S^1$,
$\varphi (\mu_1, \mu_2) =(\mu_1 \mu_2 ^q, \mu _2 ^{-p})$ is exactly
$\Gamma$.  This gives us an isomorphism $\bar{\varphi}: \bbT^2/\Gamma
\to \bbT^2$.  With this identification the image of the moment map for
the action of $\bbT^2 \simeq \bbT^2/\Gamma$ is
$$
C_{p,q} := \{ t_1 (1, 0) + t_2 (p, q) \in \RR ^2 \mid  t_1, t_2 \geq 0\}.
$$
Note that if we pick a different basis of the weight lattice of
$\bbT^2$, the moment cone $C_{p,q}$ will change by an action of an
element of $\SL (2, \ZZ)$.

We claim that we can obtain the cone $C_{p,q}$ (up to a change of
lattice basis) by cutting the image of $T^*_0 S^2$ with a half-space. Indeed
 the matrix $\left( \begin{array}{cc}
1 & 1\\
1 &2
\end{array} \right)$ maps $C_{p, q} $ onto
$$
C' _{p, q} :=
        \{ t_1 (1, 1) + t_2 (p+2q, p+ q) \in \RR ^2 \mid  t_1, t_2 \geq 0\}
$$
and
$$
C' _{p, q} = C \cap \{ \eta \in \RR^2 \mid
\langle \eta , (p+2q, -p -q)\rangle \geq 0 \},
$$
where, as above, $C$ denotes the moment image of $T^*_0 S^2$.  We
conclude that there is a Hamiltonian action of $\bbT^2$ on $L(p, q)
\times \RR$ such that the moment map image is the cut of the moment
map image of $T^*_0 S^2$ by a half-space.  It follows from
Remark~\ref{last_rem} that
$$
\left( T^*_0 S^2 \right)_{+ (p+2q, -p -q)} = L(p, q) \times \RR ,
$$
i.e., that we can obtain the symplectic cone on the lens space
$L(p,q)$ by cutting the punctured cotangent bundle of $S^2$.

More generally almost all toric symplectic cones can be obtained by
iterated cuts starting with the cotangent bundle of the standard
$n$-torus $\bbT^n$. Indeed, as remarked above, strictly convex
rational polyhedral cones in $\RR^n$ (satisfying certain integrality
conditions) classify, as moment map images, a large class of
symplectic toric cones.  Each of these polyhedral cones is the
intersection of finitely many half-spaces with primitive integral
normals. Therefore these moment map images can be obtained from $\RR^n
\smallsetminus \{0\}$ by repeated cuts by half-spaces.  Consequently
the corresponding symplectic cones can be obtained from the punctured
cotangent bundle of the standard torus $\bbT^n$ by repeated symplectic
cuts.

\end{document}